\documentclass[10pt]{amsart}

\usepackage{latexsym, amssymb, amsmath, amscd}
\usepackage{amsmath}

\input xy
\xyoption{all}

\newtheorem{theorem}{Theorem}[section]
\newtheorem*{theorem*}{Theorem}
\newtheorem{lemma}[theorem]{Lemma}
\newtheorem{remark}[theorem]{Remark}
\newtheorem{example}[theorem]{Example}
\newtheorem{corollary}[theorem]{Corollary}
\newtheorem*{corollary*}{Corollary}

\newtheorem{definition}[theorem]{Definition}

\newcommand{\adot}{\mathbf{a}_{\bullet}}
\newcommand{\bdot}{\mathbf{b}_{\bullet}}
\newcommand{\bfa}{\mathbf{a}}
\newcommand{\bfb}{\mathbf{b}}

\newcommand{\Ba}{\mathbf{B}( \mathbf{a}_{\bullet})}
\newcommand{\J}{\mathcal{J}}
\newcommand{\Ox}{\mathcal{O}_X}
\newcommand{\N}{\mathbf{N}}

\newcommand{\ord}{\textnormal{ord}}
\newcommand{\spn}{ \langle n \rangle}
\newcommand{\Spec}{\textnormal{Spec}}

\newcommand{\m}{\ensuremath{\mathfrak m}} 

\newcommand{\NN}{\ensuremath{\mathbb N}}

\newcommand{\CC}{\ensuremath{\mathbb C}}

\newcommand{\ga}{\ensuremath{\mathfrak a}}

\newcommand{\set}[1]{\ensuremath{ \left\{\, #1\, \right\} }}

\newcommand{\deq}{\ensuremath{ \stackrel{\textrm{def}}{=}}}

\begin{document}

\title{A Brian\c con--Skoda type theorem for graded systems of ideals}

\author{Alex K\"uronya, Alexandre Wolfe}

\begin{abstract}
We establish a generalization of the Brian\c con--Skoda theorem about  integral closures 
of ideals for graded systems of ideals satisfying a certain geometric condition.
\end{abstract}

\maketitle

\bigskip

\section{Introduction}

A celebrated result of Brian\c con and Skoda \cite{BS}
states that for every coherent 
sheaf of ideals $\ga \subseteq \Ox$ on a smooth irreducible complex variety, 
a certain fixed power of its integral closure is actually contained in $\ga$. 
We establish an analogue of this statement for stable graded systems of 
ideals.   

The Brian\c con--Skoda theorem, originally proved via $L^2$-methods in complex 
ana\-lysis, has sparkled a great deal of interest among algebraists. Algebraic 
proofs were given by Lipman and Teissier for ideals in a regular local 
ring \cite{LipmanTeissier}, Hochster and Huneke \cite{HochsterHuneke} 
using tight closure theory (one of the first major victories of this method), 
and Ein--Lazarsfeld \cite{GeomNSS, PAG}, via multiplier ideals. 
Here we will follow this latter path, and prove our result by means of 
multiplier ideal theory. For a thorough treatment of this topic, the reader 
is invited to  consult  \cite[Chapters 9--11]{PAG}.

Let $X$ be a smooth irreducible complex variety.
A graded system of ideals on $X$ 
is a sequence ideal sheaves
${\ga}_m\subseteq \Ox$, for which $\bfa_m \bfa_n \subseteq \bfa_{m+n}$ for 
all natural numbers $m,n\geq 1$. Experience shows (see \cite{AWthesis})
that unlike the powers of a fixed ideal,  graded systems of ideals
 in general behave  very pathologically.  In an attempt to remedy 
this situation, we introduce the concept of stability, and show 
that stable graded systems of ideals act in a more controlled manner.

We define a stable graded system of ideals to be one, where all 
elements differ from 
$\Ox$, and for every irreducible subvariety $Z$ of $X$, which occurs  
in the support of $\bfa_m$ for infinitely many $m$, the asymptotic 
order of vanishing of $\adot$ along $Z$ is positive. This notion is in a 
way an algebraic counterpart of the notion of a stable divisor (for a 
precise definition of this notion the reader is invited 
to consult \cite{AIBL}); the graded 
system of ideals formed by the base ideal sheaves of the multiples of a stable
Cartier divisor is certainly stable.  

A distinguishing property of stable graded systems of ideals is  that their cosupports
---  and the cosupports of the corresponding sequence of asymptotic multiplier ideals ---
stabilize to an algebraic set. 

Our main result is the following. 

\begin{theorem*}
Let $\adot$ be a stable graded system of ideals.  Then there exists an integer $C$  so that for all $n \gg 0$, 
$$
\J ( Cn\cdot \adot) \subseteq \bfa_n.
$$
\end{theorem*}

As asymptotic multiplier ideals are integrally closed, the fact that 
$\overline{\bfa}_m\subseteq \J (m\cdot \bfa.)$ then implies  the following 
generalization of the aforementioned result of Brian\c con and Skoda. 

\begin{corollary*}
Let $\adot$ be a stable graded system of ideals. Then there exists a positive integer 
$C$, such that for all $n\gg 0$
\[
\overline{\bfa}_{Cn}\subseteq  \bfa_n\ .
\]
\end{corollary*}

Note that as a by-product, our results give a relatively elementary treatment
(modulo standard properties of asymptotic multiplier ideals) of a weak version 
of the original  Brian\c con--Skoda theorem with the exception of finitely many initial terms. 

The organization of the paper goes as follows. The definition and basic properties of 
stable graded systems of ideals are found in Section 2, while Section 3 hosts the generalized
Brian\c con--Skoda theorem along with its proof and some illustrating examples. 

We would like to thank Manuel Blickle, Lawrence Ein, and Robert Lazarsfeld for helpful discussions.

\section{Stable Graded Systems of Ideals}

Throughout this paper we will work on a smooth  irreducible complex variety $X$ of dimension $d$.

\begin{definition}
Let $\adot = \{ \bfa_n \}_{n \in \N}$ be a  sequence of ideal sheaves on  $X$.  We say that $\adot$ is a {\em graded system of ideals} 
if for all $m, n \in \N$, 
$$
\bfa_m \bfa_n \subseteq \bfa_{m+n}.
$$
\end{definition}

Setting $\bfa_0 = \Ox$, we note that the definition is equivalent to asking that  $\oplus_{n \geq 0} \bfa_n$ have a natural structure of a graded algebra of sheaves on $X$. This is called the Rees algebra associated to the graded system $\adot$.
Notable examples of graded systems of ideals  are the system of symbolic powers of a prime ideal,  
or the system formed by the  base ideals of the tensor powers of a line bundle on a projective variety.

If $\adot$ is a graded system of ideals, the set 
\[
\{ n \in \N : \bfa_n \neq 0 \}
\]
is a subsemigroup $S(\adot)$ of $\N$.  It is well-known that this forces all 
sufficiently large elements of $S(\adot)$ to be multiples of 
$e( \adot) \deq \gcd S( \adot)$.  We call $S(\adot)$ the semigroup of 
$\adot$, and $e(\adot)$ the exponent of $\adot$.  In the following, we shall 
deal exclusively with graded systems of ideals with exponent $e(\adot) = 1$, 
so that $\bfa_n \neq 0$ for all $n \gg 0$.

Let $Z \subset X$ be a (reduced and irreducible) subvariety, and $\bfa$ an ideal sheaf on $X$.  Then
the order of vanishing of the ideal $\bfa $ along the subvariety $Z$ is 
\[
\ord_Z (\bfa) \deq \max \{ n \in \N : \bfa \subseteq I_Z^{\spn} \} \ ,
\]
where $I_Z^{\spn}$ denotes the $n$th symbolic power of the ideal $I_Z$.

If $Z$ is a point $p \in X$, then $\ord_p (\bfa) = n$ simply says that all partial derivatives of order at most $n-1$ of functions in $\bfa$ vanish at $p$.  In fact, the case $Z = p$ is already in some sense the general case, since
\[
\ord_Z (\bfa) = \min \{ \ord_p (\bfa) : p \in Z \} \ ,
\]
and this minimum is attained at general points $p \in Z$.  Along similar lines, one can see that it suffices to test order of vanishing generically: $\ord_Z (\bfa) = n$ if and only if $\ord_p (\bfa) = n$
for general $p \in Z$, $\ord_p (\bfa) = n$.

Recall that \cite{AIBL} shows that one may attach an asymptotic order of vanishing to a graded system of ideals $\adot$ via the following formula:

\begin{definition}
Let $\adot$ be a graded system of ideals.  Then 
$$
\ord_Z (\adot) := \lim_{n \rightarrow \infty} \frac{ \ord_Z ( \bfa_n )}{n}.
$$
\end{definition}
\noindent
More precisely, it is established in \cite{AIBL} that the limit in question
always exists. 

We now arrive at the key technical definition.

\begin{definition}
We say that a  graded system $\adot$ is {\em stable} if 
\begin{enumerate}
\item $e(\adot) = 1$ \label{exp1}
\item for all subvarieties $Z$ that occur as components of $V(\bfa_n)$ for infinitely many $n \in \N$, $\ord_Z (\adot) > 0$. \label{positiveorder}
\end{enumerate}
\end{definition}

Stable graded systems of ideals are ubiquitious. Ordinary or symbolic powers of a prime ideal certainly form stable graded systems, and various nontrivial examples of them occur naturally in geometric contexts.

\begin{remark}{\rm 
Let $X$ be a smooth projective variety, $D$ a stable Cartier divisor 
(in the sense of \cite{AIBL}) on $X$. Then the graded system of ideals 
${\mathfrak b}_m(D)= {\mathfrak b}(mD)$ is stable 
(see \cite[Proposition 3.8.]{AIBL}).
}
\end{remark}

\begin{example}{\rm 
The following is a classical example of a non-stable graded system of ideals 
which goes back to Zariski. In \cite{Zariski} he shows the existence of  
 a smooth projective surface $X$, an irreducible 
curve $C\subseteq X$, and a divisor $L$ on $X$, such that for every $m\geq 1$, the curve  $C$ 
is in the base locus of $mL$, but $mL-C$ is base point free. In this case,  the base ideal $\mathfrak{b}(mL)$   associated to $mL$ vanishes to order exactly one along $C$, hence the graded system
 consisting of the base ideals $\mathfrak{b}(mL)$ is not stable.
}
\end{example}

A key fact about stable graded systems is that their cosupports stabilize 
to an algebraic set, the asymptotic cosupport of $\adot$, which we will 
denote $\Ba$. Later we show that $\Ba$ also equals the asymptotic cosupport
of the corresponding sequence of asymptotic multiplier ideals.

\begin{lemma}
Let $\adot$ be a stable graded system of ideals.  Then there exists an algebraic set $\Ba \subset X$ so that for all $n \gg 0$, $V(\bfa_n) = \Ba$.
\end{lemma}

\begin{proof}  
In fact, we will show that 
$$
\Ba = \bigcap_{n \in \N} V ( \bfa_n).
$$
Since $X$ is a noetherian topological space, there is some $N \in \N$ so 
that $\Ba = \cap_{i=1}^N V( \bfa_i)$.  In particular, 
$\Ba \subseteq V(\bfa_n)$ for all $n \geq N$.  Stability will guarantee 
the reverse inclusion.

We first show that the subsequence $n_k = k \cdot N!$ satisfies $V(\bfa_{n_k}) = \Ba$ for all $k$.  By definition of graded systems, for all $i$, $1 \leq i \leq N$, one has $\bfa_i^{N!/i} \subseteq \bfa_{N!}$, so $V(\bfa_{N!}) \subseteq V(\bfa_i)$.  Thus 
$$
V(\bfa_{N!}) \subseteq \bigcap_{i=1}^N V(\bfa_i) = \Ba.
$$
So for all $k \in \N$, one has
\[
V(\bfa_{k \cdot N!}) \subseteq V(\bfa_{N!}^k) = \Ba \ ,
\]
as desired.

Now, let $p \in \NN$ be a prime greater than $N$.  Then $p$ and $N!$ generate 
a subsemigroup of $S(\adot)$ under addition which includes all $n \gg 0$.  
Hence one has $n = n_1 p + n_2 N!$ for all $n\gg 0$, and 
\[
V(\bfa_n) \subseteq V(\bfa_{p}^{n_1}) \cup V(\bfa_{N!}^{n_2}) = 
V(\bfa_p) \cup V(\bfa_{N!})\ .
\]

We now claim that there exists $n_0\in {\NN}$ depending on $\adot$
and $p$ such that if $n\geq n_0$ then in fact 
\[
V(\bfa_n) \subset V \cup V(\bfa_{N!})\ ,
\]
where $V \subset V(\bfa_p)$ and dim$(V) < $ dim$(V(\bfa_p))$.  

Granting the claim, we may then repeat the same argument with $V(\bfa_p)$ 
replaced by $V$; since $V(\bfa_p)$ is itself a Noetherian topological space, 
after finitely many iterations of this process we obtain that 
$V(\mathbf{a}_n) = \Ba$ for $n \gg 0$.

To prove the claim, let $Z$ be a top-dimensional component of $V(\bfa_p)$ which does not occur in $V(\bfa_{N!})$.  Then $Z$ does not appear in $V(\bfa_{k \cdot N!})$ for all $k$, so calculating along that subsequence we see that $\ord_Z (\adot) = 0$.  By definition of stability, we see that $Z$ can occur as a component of $V(\mathbf{a_n})$ for at most finitely many $n$.  So after finitely many steps $V(\mathbf{a_n})$ is supported along some proper algebraic subset of $Z$.  Repeating this argument for each top-dimensional component of $V(\bfa_p)$, the claim follows. 
\end{proof}

Recall that  the level-$n$ asymptotic multiplier ideal $\J (n \adot)$ is 
defined to be the maximal element of the set of multiplier ideals 
\[
\set{ \J( \frac{n}{p}\cdot {\bfa}_p )}\ .
\]
For a proof of the fact that this maximum exists and other properties 
of asymptotic multiplier ideals, the reader is referred to 
\cite[Chapter 11]{PAG}. 

One thus obtains a sequence of ideals 
$$
\bfb_n := \J (n \adot),
$$
indexed by the natural numbers. It is not graded; however, the subadditivity
theorem for asymptotic multiplier ideals shows that the ideals $\bfb_n$ form  a reverse graded
system (in the sense of \cite{Mustata}):
\[
\bfb_{n+m} \subseteq \bfb_n \cdot \bfb_m \ ,
\]
which we will denote by $\bdot$.
In \cite{Mustata} 
it is shown that the following limit defining $\ord_Z (\bdot)$ exists:
$$
\ord_Z(\bdot) := \lim_{n \rightarrow \infty} \frac{\ord_Z (\bfb_n)}{n} \ .
$$
Analogously to the case of stable graded systems of ideals, we now verify that the cosupports of an asymptotic  multiplier ideal sequence corresponding to a stable graded system also stabilize.

\begin{lemma}
Let $\adot$ be a stable graded system of ideals with asymptotic cosupport $\Ba$.  Then for all $n \gg 0$, $V( \J (n \adot)) = \Ba$.
\end{lemma}

\begin{proof}  
By the definition of multiplier ideals one has that for all 
$n \in \N$, $\bfa_n \subset \J (n \adot)$.  Therefore 
$V( \J (n \adot)) \subseteq \Ba$ for $n\gg 0$, and 
\[
\ord_Z (\bdot) \leq \ord_Z(\adot)
\]
for  all subvarieties $Z$ of $X$.  We claim that in fact 
\[
\ord_Z (\bdot) = \ord_Z (\adot)
\]  
for every  subvariety $Z$ of $X$.
This will show in particular that for $n \gg 0$, $\J ( n \adot)$ is co-supported along every irreducible component of $\Ba$ and hence $V( \J ( n \adot)) = \Ba$.  To establish the claim, the key will be to compare orders of vanishing on $X$ and on suitable log resolutions of the ideals in the graded system.

Fix a subvariety $Z\subseteq X$, and 
let $X_0 = X$ and for all $i\geq 1$, pick 
\[
\mu_i: X_i \rightarrow X 
\]
to  be log resolutions of $\bfa_i$ such that 
\begin{enumerate}
\item ${\mu}_i$ is a log resolution of $Z$,
\item for all $j<i$, ${\mu}_i$ factors through ${\mu}_j$.
\end{enumerate}


As  the $\mu_i$'s are log resolutions of $Z$ for all $i$,  $\mu_i^{-1}(Z)$ is
a SNC divisor on $X_i$ for all $i$, and  there is a unique prime divisor $E_i \subseteq X_i$ dominating the smooth locus of $Z$.Write $b_i$ for the coefficient of $E_i$ in the relative canonical divisor $K_{X_i/X}$, and $a_i$ for the coefficient of $E_i$ in $F_i = \mu_i^{*} \bfa_i$.

Observe firstly that $b_i = b_1$ for all $i \geq 0$, as $\mu_1^{-1}(Z)$ is already a SNC divisor on $X_1$.  Then for some $p = p(n)$ which we may take to be arbitrarily large,
\begin{eqnarray*}
\ord_Z (\J (n \adot)) &\geq& \textnormal{ord}_{E_p} (K_{X_p/X} - [ \frac{n}{p} F_i]) \\
&=& -b_1 + [n \frac{a_p}{p}] \\
&=& -b_1 + [n \frac{ \ord_Z (\bfa_p)}{p}].
\end{eqnarray*}
Fix a sequence $p(n)$ as above, with the additional property that 
$p(n)\rightarrow\infty$ as $n\rightarrow\infty$.
Estimating the asymptotic order of $\bdot$ gives
\begin{eqnarray*}
\ord_Z (\bdot) &\geq& \lim_{n \rightarrow \infty}
 \frac{-b_1 + [n \tfrac{ \ord_Z (\bfa_{p(n)})}{p(n)}] }{n} \\
& = & \lim_{n\rightarrow\infty}{ \frac{n\cdot \tfrac{\textnormal{ord}_{Z}( 
\bfa_{p(n)})}{p(n)} - \left\{ n \tfrac{ \ord_Z (\bfa_{p(n)})}{p(n)} \right\}}
{n}} \\
& = &  \lim_{n\rightarrow\infty}{ \frac{n\cdot \tfrac{\textnormal{ord}_{Z}( 
\bfa_{p(n)})}{p(n)}}{n}}  \\
 &=& \ord_Z (\adot)\ , 
\end{eqnarray*}
as we wanted.
\end{proof}


\section{A generalized statement of Brian\c con-Skoda type}

In analogy with the aformentioned result of Brian\c con and Skoda, 
we are interested in finding functions $f: \N \rightarrow \N$ so that 
for all $n \gg 0$,
$$
\overline{\bfa}_{f(n)} \subseteq \bfa_n\ .
$$
Our plan is to establish the corresponding statement, ie.
\[
\J (\bfa_{f(n)}) \subseteq \bfa_n
\]
if $n\gg 0$, for asymptotic multiplier ideals, and appeal to their integral closure. 

We note first of all that such functions may not exist if the graded system $\adot$ is not stable.

\begin{remark}{\rm  
Take $X = {\CC}^2$, $\ell \subset X$ a line and $p \in \ell$ a point.  Let $\adot$ be a graded system generated by $\bfa_1$ and $\bfa_2$, where $\bfa_1 = I_{\ell}$ and $\bfa_2 = I_p$.  Then
\[
\bfa_n =  \begin{cases}  I_p^{n/2} & \textnormal{ if } n \textnormal{ is even} \\
                  I_p^{(n-1)/2}I_{\ell} & \textnormal{ if } n \textnormal{ is odd}. 
\end{cases}
\]
For any $k$, we may always calculate the level-$k$ asymptotic multiplier ideal at an ideal $\bfa_{p(k)}$ where $p(k)$ is even, which shows that $\J (k \adot)$ is cosupported at $p$.  Therefore, if $n$ is any odd integer, for all $k \in \NN$, $\J (k \adot) \not\subseteq \bfa_n$.

We remark that, as we have seen, stability of $\adot$ rules out the troublesome appearances of $I_{\ell}$ to bounded order in $\mathbf{a}_n$ for infinitely many $n$.  $\Box$
}
\end{remark}

On the other hand, in the case of stable graded systems of ideals, we may find functions $f$ which have linear growth and satisfy the desired containment.

\begin{remark}{\rm 
Consider first the case of finitely generated graded systems, ie. graded systems of ideals for 
which the associated Rees algebra $\bigoplus_{m\geq 0}{\bfa_m}$ is finitely generated. Then by 
Proposition 2.4.27. in \cite{PAG}, there exists an integer $N>0$ such that for all $m\geq 1$
one has 
\[
{\bfa}_{N}^{m} = {\bfa}_{Nm}\ .
\]
If  $\adot$ is stable then 
\[
\J ((Nm+D)\cdot \adot ) \subseteq \bfa_m
\]
for a suitable positive constant $D$. Note that finitely generated graded systems are not 
necessarily stable.
}
\end{remark}

\begin{theorem} \label{skoda}
Let $\adot$ be a stable graded system of ideals.  Then there exist positive constants $C$ and $D$ so that for all $n \gg 0$, 
\[
\J ( (\lceil Cn  + D\rceil ) \adot) \subseteq \bfa_n.
\]
In particular, 
\[
\J ( C'n\cdot \adot ) \subseteq \bfa_n
\]
for a suitable positive integer $C'$.
\end{theorem}
\begin{proof}  
Fix two relatively prime integers $n_1 < n_2$, so large that 
$V(\bfa_i) = V( \J( n_1 \adot)) = \Ba$ for $i=1,2$.   The semigroup generated 
by the $n_i$  includes all $n \gg 0$.  By Hilbert's Nullstellensatz, there 
 exist integers $\nu_i$ so that for $i=1,2$
$$
I_{\Ba}^{\nu_i} \subseteq \bfa_{n_i}.
$$
We are looking for a function $f:\NN\rightarrow \NN$ such that  $f$ satisfies 
$\J (f(n) \adot) \subset \bfa_n$ for all $n$ in the semigroup generated by the
$n_i$.  Observe that 
\[
\J ( f(n) \adot) \subseteq \J (n_1 \adot)^{ [ f(n)/n_1]} 
\subseteq  I_{\Ba}^{[ f(n)/n_1 ]}.
\]
Write $n = m_1 n_1 + m_2 n_2$.  Then if $[f(n)/n_1] \geq m_1 \nu_1 + m_2 \nu_2,$
\[
\J (f(n) \adot) \subseteq I_{\Ba}^{m_1 \nu_1 + m_2 \nu_2} 
 \subseteq \bfa_1^{m_1} \bfa_2^{m_2} 
 \subseteq \bfa_n.
\]
Since in any case the $m_i \leq \tfrac{n}{n_i}$, one checks that $f(n) =
\lceil (\nu_1 + \tfrac{\nu_2 n_1}{n_2})n \rceil + n_1$ satisfies the
containment. 
\end{proof}

\begin{corollary}
Let $\adot$ be a stable graded system of ideals. Then there exists a positive integer $C$, 
such that for all $n\gg 0$, 
\[
\overline{\bfa}_{Cn}\subseteq  \bfa_n\ .
\].
\end{corollary}

It is natural to ask what numbers $C$ and $D$ may occur in Theorem \ref{skoda}.  In many situations
of geometric interest, for instance,  $C=1$ will do.  
However, simple examples show that this is not  necessarily the case in
general. 

\begin{example}{\rm
 Let $X = \Spec (\CC [x] )$ and $I = (x)$.  One checks that $\bfa_n = I^{n +
 \lceil \sqrt{n} \rceil}$ is a graded system of ideals.  However, $\J ( (n+D)
 \adot) = I^{n+D} \not\subseteq \bfa_n$ for all $n \gg 0$. 
}
\end{example}

\begin{example}{\rm
Here we show that in certain cases the set of possible candidates for the 
constant $C$ is bounded away from $1$.  
Let $\m = (x,y)\subseteq {\CC}^2$ and let $\epsilon > 0$ be arbitrary. Set
\[
\bfa_n = {\m}^{\lfloor\epsilon n \rfloor} (x^n, y)\ .
\]

Since $y \in ( {\bfa}_n : {\m}^{\lfloor\epsilon n\rfloor})$ 
for all $n$, the multiplier ideals of the graded
system $\adot$  and ${\m}^{\lfloor\epsilon n\rfloor}$ are the same. 
So if  $\J( (Cn+D) \adot) \subseteq {\bfa}_n$, one also has
\[
{\m}^{ \lfloor\epsilon (Cn+D)\rfloor } \subseteq {\bfa}_n\ .
\] 
However, compare the powers of $x$ that occur in both sides of the 
containment directly above. On the right-hand side, $x^{n + \lfloor\epsilon n\rfloor}$ 
is the lowest power. So   $\lfloor\epsilon (Cn+D)\rfloor \geq n + \lfloor\epsilon n\rfloor $,
and taking a limit as $n\rightarrow  \infty$, one has $C >1 + \tfrac{1}{\epsilon}$.
}
\end{example}

\begin{remark}{\rm 
Here is  what happens for sums and products of graded systems of ideals.
Let $\adot$, $\bdot$ be two stable graded systems of ideals on $X$, $C_A,C_B$ two positive 
integers which satisfy the conclusion of the theorem. If $\adot\cdot\bdot $ and 
$\adot + \bdot $ denote the product and sum of the two graded systems (in the sense of 
\cite[2.4.24.]{PAG}), then for large enough $n$, one has 
\begin{eqnarray*}
\J ( (C_A+C_B)n\cdot (\adot\cdot\bdot)) & \subseteq & (\adot\cdot\bdot)_n, \\
& \textrm{ and } & \\
\J ( (C_A+C_B)n\cdot (\adot + \bdot)) & \subseteq & (\adot + \bdot)_n\ .
\end{eqnarray*}
These results rest upon the subadditivity theorem for multiplier ideals and Mustat\c t\v a's summation theorem
 \cite{Mustata2},  respectively.
 } 
\end{remark}

\bigskip

{\small\sc
\noindent
Alex K\"uronya \\
Universit\"at Duisburg-Essen, Campus Essen, Fachbereich 6 Mathematik \\
D-45117 Essen, Germany \\
\ \\ 
Budapesti M\H{u}szaki \'es Gazdas\'agtudom\'anyi Egyetem, Matematikai Int\'ezet\\
H-1521 Budapest, P.~O.~Box 91 \\
Hungary \\
{\em email address:\ }{\tt alex.kuronya@math.bme.hu} \\ 
\ \\
Alexandre Wolfe \\
Department of Mathematics, University of Michigan \\
Ann Arbor, MI 48109-1109\\
USA \\
{\em email address:\ {\tt alexwolfe@hotmail.com} \\

}


\begin{thebibliography}{99}

\bibitem{BS}
   Skoda, Henri, Brian\c con, Jo\"el. 
   Sur la cl\^oture int\'egrale d'un id\'eal de germes de fonctions 
   holomorphes en un point de $C\sp{n}$.  
   C. R. Acad. Sci. Paris S\'er. A 278 (1974), 949--951.

\bibitem{GeomNSS}
   Ein, Lawrence,  Lazarsfeld, Robert. 
   A geometric effective Nullstellensatz. 
   Invent. Math. 137 (1999), no. 2, 427--448.

\bibitem{AIBL}
   Ein,~Lawrence, Lazarsfeld, ~Robert, Musta\c{t}\v{a},~Mircea,
   Nakamaye,~Michael, Popa,~Mihnea. 
   Asymptotic invariants of base loci.
   Preprint, \texttt{math.AG/0308116}.

\bibitem{HochsterHuneke}
   Hochster, Melvin; Huneke, Craig. 
   Tight closure, invariant theory, and the Brian\c{c}on-Skoda theorem. 
   J. Amer. Math. Soc. 3 (1990), no. 1, 31--116.

\bibitem{PAG}
  Lazarsfeld,~Robert.  
  Positivity in Algebraic Geometry I.-II. Ergebnisse der Mathematik und ihrer 
  Grenzgebiete, Vols. 48-49., Springer Verlag, Berlin, 2004.

\bibitem{LipmanTeissier}
  Lipman, Joseph, Teissier, Bernard.
  Pseudorational local rings and a theorem of 
  Brian\c{c}on-Skoda about integral closures of ideals. 
  Michigan Math. J. 28 (1981), no. 1., 97--116.

\bibitem{Mustata}
  Musta\c t\v a, Mircea. 
  On multiplicities of graded sequences of ideals. 
  J. Algebra 256 (2002), no. 1., 229--249.
 
\bibitem{Mustata2}
  Musta\c t\v a, Mircea.
  The multiplier ideals of a sum of ideals.  
  Trans. Amer. Math. Soc. 354 (2002), no. 1, 205--217.  
   
  
\bibitem{AWthesis}
  Wolfe, Alexandre.
  Ph.D. thesis,
  University of Michigan, 2005.  
  
\bibitem{Zariski}
  Zariski, Oscar.
  The theorem of Riemann--Roch for high multiples of an effective divisors on an 
  algebraic surface.
  The Annals of Math. 2nd Ser., 76 (1962), no. 3., 560--615.   

\end{thebibliography}
\end{document}